\documentclass[12pt]{article}
\usepackage{mathrsfs}
\usepackage[dvips]{epsfig}
\usepackage{amsmath}
\usepackage{amsfonts}
\usepackage{amssymb}
\usepackage{graphicx}

\begin{document}

\title{A new compact class of open sets under Hausdorff distance and shape optimization }
\author{Donghui Yang \\
School of Mathematics and Statistics, Central China Normal University,\\ 430079, Wuhan, People's Republic of China
\\ \textsf{donghyang@mail.ccnu.edu.cn}\\}
\date{}
\maketitle

\begin{abstract}
In this paper we obtain a  new class of open sets, and we prove the
class is compact under the Hausdorff distance, then we prove the existence
of solutions of some shape optimization for elliptic equations.
\end{abstract}

\thispagestyle{empty}


\baselineskip=17.5pt \parskip=3pt

\section{ Introduction}

\hskip600pt \footnote{The author was supported by NSFC 10901069.}

There are many papers concerning existence theory for shape optimization problems.
There are several types of results: using regularity assumptions for the boundary
of the unknown open sets (see \cite{BU}, \cite{ch1}, \cite{p1}, \cite{t1}), using certain capacitary
constraints (see \cite{B4}, \cite{h1}, \cite{s1}) or using the notion of a generalization perimeter
and constraints or penalty terms constructed with it (see \cite{B2}, \cite{B3}, \cite{B4}, \cite{h1}).
In general, for the shape optimization problems, we must give a class of open sets and
prove the class is complete under the Hausdorff distance in the first place, and then
we should prove that the shape optimization problems have at least one solutions.

In this work, we give a  new class of open sets, and we prove the
class is compact under the Hausdorff distance, then we prove the existence
of solutions of some shape optimization for elliptic equations.

Let $B^*\subset\mathbb{R}^k$ be a bounded domain. On the space $H^1_0(B^*)$ we consider
the norm $\displaystyle\|u\|_{H_0^1(B^*)}={\Big (}\int_{B^*}|\nabla u|^2dx{\Big )}^{\frac{1}{2}}$. Let be given a
smooth symmetric matrix $A\in M_{k\times k}(C^1(\overline{B^*}), A=A^T$
(where $A^T$ represent the transformation of the matrix $A$), and
\begin{eqnarray}\label{0}
\alpha |\xi|^2\leq \langle A.\xi,\xi\rangle,
\end{eqnarray}
where $0<\alpha$ is a constant. We define the
associated operator $\mathscr{A}: H^1_0(B^*)\rightarrow H^{-1}(B^*)$:
\begin{eqnarray}\label{1}
\mathscr{A}=\mbox{div}(A.\nabla).
\end{eqnarray}

For any open set $\Omega\in\mathscr{C}_{M,R}$ ($\mathscr{C}_{M,R}$ will define later),
we consider the Dirichlet problem in $\Omega$:
\begin{eqnarray}\label{2}
u\in H_0^1(\Omega),\quad -\mathscr{A} u_\Omega=f
\end{eqnarray}
in the variational sense, i.e.,
\begin{eqnarray}\label{3}
\int_\Omega \langle A.\nabla u_\Omega,\nabla\phi\rangle dx
=\langle f|_\Omega,\phi\rangle_{H^{-1}(\Omega)\times H_0^1(\Omega)}\quad \forall \phi\in C_0^\infty(\Omega)
\end{eqnarray}
with $f\in H^{-1}(B^*)$ and $f|_\Omega$ denoting the restriction of the distribution $f$
to the open set $\Omega$.

Because $H_0^1(\Omega)=cl_{H_0^1(B^*)}(C_0^\infty(\Omega))$, then (\ref{3}) has a unique solution
$u_\Omega\in H_0^1(\Omega)$, which we can extend with zero on $B^*-\Omega$, to $u_\Omega^0$, and
$u_\Omega^0\in H_0^1(B^*)$, $\|u_\Omega^0\|_{H_0^1(B^*)}=\|u_\Omega\|_{H_0^1(\Omega)}$.

When we consider the solution of (\ref{2}), we will implicitly take its extension $u_\Omega^0$. We
shall study the following shape optimization problem
$$\inf_{\Omega\in\mathscr{C}_{M,R}} J(\Omega)
\equiv\inf_{\Omega\in\mathscr{C}_{M,R}}{\Big\{} \frac{1}{2}\int_{B^*}|u_\Omega-g|^2dx{\Big \}}.\eqno{(P)}$$

In the following, we give the definition of $\mathscr{C}_{M,R}$.

We call a open set $\Omega\subset\mathbb{R}^k$ with property ($C_M$),
if for any $x, y\in \Omega$, there exists a connected compact set
$K$ with $x, y\in K$, such that $K\subset\Omega$ and
$\displaystyle\bigcup_{z\in K}B(z, \frac{d^*}{M})\subset\Omega$,
here $d^*=\min\{\mbox{dist}(x,\partial\Omega), \mbox{dist}(y,
\partial\Omega)\}$ and $M>1$ is a given constant. Here and throughout this paper we
denote $B(z, r)$ an open ball with center $z$ and radius $r$ .

For given $R>0$ small enough, we define
\begin{eqnarray*}
\mathscr{C}_{M,R}=\{\Omega\subset\mathbb{R}^k;\ B(x_\Omega,
R)\subset\Omega\subset B^*, \Omega \ \mbox{is a open set and satisfy
the property}\ (C_M)\},
\end{eqnarray*}
where $B^*$ is a bounded domain and $x_\Omega$ is some point in
$\Omega$.

The topology on $\mathscr{C}_{M,R}$ is induced by the
Hausdorff distance between the complementary sets, i.e., for
any $\Omega_1, \Omega_2\in \mathscr{C}_{M,R}$,
\begin{eqnarray}\label{1.5}\displaystyle
\rho(\Omega_1, \Omega_2)=\mbox{max}\{\sup_{x\in \overline
{B^*}\setminus\Omega_1} \mbox{dist}(x, \overline
{B^*}\setminus\Omega_2), \sup_{y\in \overline
{B^*}\setminus\Omega_2} \mbox{dist}(\overline
{B^*}\setminus\Omega_1, y)\}, \;\;\;
\end{eqnarray}
where $\mbox{dist}(\cdot,\cdot)$ denotes the Euclidean metric in
$\mathbb{R}^N$. We denote by Hlim, the limit in the sense of
(\ref{1.5}).

In this work, we obtain the following main result about the family $
\mathscr{C}_{M,\rho}$:

\vskip 1mm

{\bf Theorem 2.1} If $\{\Omega_m\}_{m=1}^\infty\subset
\mathscr{C}_{M,R}$, then there exists a subsequence
$\{\Omega_{m_k}\}_{k=1}^\infty$ of
$\displaystyle\{\Omega_m\}_{m=1}^\infty$ such that
$$\mbox{H}\!\!\!
\lim_{k\rightarrow \infty}\Omega_{m_k}=\Omega\ \ \mbox{and}\ \
\Omega\in \mathscr{C}_{M,R}.$$ i.e., $(\mathscr{C}_{M,R}, \rho)$ is
a compact metric space.

\vskip 1mm

We note that there are many compact classes of open sets under Hausdorff
distance have been found (See \cite{B2}, \cite{B3},
\cite{BU}, \cite{ch1}, \cite{ch2}, \cite{MM}, \cite{s1}, \cite{t1},
\cite{wy} etc.). We also can see more cases in the books \cite{B4},
\cite{h1}, \cite{p1} etc. But the class of the open sets $\mathscr{C}_{M,R}$
 has never appeared in any other place.

By Theorem 2.1, we obtain the existence of the optimal solutions for
problem $(P)$:

\vskip 1mm

{\bf Theorem 3.1}\quad The shape optimization problem $(P)$ has at
least one solution.

\section {Compactness of $\mathscr{C}_{M,R}$}

\quad\ \ We shall use the following notations:
$$\displaystyle\delta(K_1,K_2)=\mbox{max}\{\sup_{x\in K_1}\mbox{dist}(x, K_2), \sup_{y\in
K_2}\mbox{dist}(K_1,y)\},$$ where $K_1$ and $K_2$ are compact subsets in
$\mathbb{R}^k$. Following from the definitions of $\rho, \delta$, we
obtain that $\rho(\Omega_1,\Omega_2)=\delta(\overline {B^*}\setminus
\Omega_1, \overline {B^*}\setminus \Omega_2)$ for any open sets
$\Omega_1, \Omega_2\subset B^*$, hence we also call $\delta$ the
Hausdorff distance.

The following results were given and proved in \cite{B4}, \cite{h1},
\cite{p1}. Which will be used in this paper.

\vskip 1mm

{\bf Lemma 2.1 } Let $A, A_n,n=1,2,\cdots$, be compact subsets in
$\mathbb{R}^k$ such that $\delta(A_n,A)\rightarrow 0$, then $A$ is
the set of all accumulation points of the sequences $\{x_n\}$ such
that $x_n\in A_n$ for each $n$.

\vskip 1mm

{\bf Remark 2.1 } Following from Lemma 2.1 and the definition of
$\delta$, we obtain that for any $\epsilon>0$ there exists
$N(\epsilon)>0$ such that for all $m\geq N(\epsilon)$ we have
$\displaystyle A\subset\bigcup_{x\in A_m}B(x, \epsilon)$ and
$\displaystyle A_m\subset\bigcup_{x\in A}B(x, \epsilon)$.

\vskip 1mm

{\bf Lemma 2.2 } Let $A,\tilde A, A_n, \tilde A_n,n=1,2,\cdots$, be
compact subsets in $\mathbb{R}^k$ such that
$\delta(A_n,A)\rightarrow 0$ and $\delta(\tilde A_n,\tilde
A)\rightarrow 0$. Suppose that $A_n\subset \tilde A_n$ for each $n$,
Then $A\subset \tilde A$.

\vskip 1mm

{\bf Lemma 2.3 } ($\Gamma-$property for $\mathscr{C}_{M,R}$)\quad
Assume that $\{\Omega_n\}_{n=1}^\infty\subset \mathscr{C}_{M,R},
\Omega_0\in \mathscr{C}_{M,R}$ and $\Omega_0=\mbox{Hlim} \Omega_n$.
Then for each open subset $K$ satisfying $\overline K\subset
\Omega_0$, there exists a positive integer $n_K$ (depending on $K$)
such that $\overline K\subset \Omega_n$ for all $n\geq n_K$.

\vskip 1mm

{\bf Lemma 2.4 } If $\Omega_n\subset B^*$, $n\in\mathbb{N}$, are open bounded sets, there exists
$\Omega\subset B^*$, open, such that $\Omega=\mbox{Hlim}\Omega_n$, on a subsequence.
In particular, let $\mathscr{O}=\{C\subset B^*;\ C\not=\emptyset,\
C\ \mbox{is compact}\}$. Then $(\mathscr{O}, \delta)$ is a compact
metric space.

\vskip 1mm

The next Lemma 2.5 is clearly, but we show it in the following.

\vskip 1mm

{\bf Lemma 2.5}\quad Let $\{x_n\}\subset \mathbb{R}^k$ and
 $\{r_n\}\subset \mathbb{R}$ be such that $r_n\geq r_0>0$ and
 $\mbox{Hlim}B(x_n,r_n)=D\subset \mathbb{R}^k$. Then there exists $x_0\in
 \mathbb{R}^k$ such that $B(x_0,r_0)\subset D$ and
 $x_n\rightarrow x_0 (n\rightarrow \infty)$.

\vskip 0.1cm

{\bf Proof}\quad Since $B(x_n,r_n)\subset B^*$ and $B^*$ is a
bounded subset in $\mathbb{R}^k$, there exists a subsequence of
$\{x_n\}$, still denoted by itself, such that
\begin{eqnarray}\label{2.9}
x_n\rightarrow x_0
\end{eqnarray}
for some $x_0\in \overline {B^*}$.

We claim that
\begin{eqnarray}\label{2.10}
B(x_0,r_0)\subset D.
\end{eqnarray}
By contradiction, we assume that there did exist $y_0\in B(x_0,r_0)$
and $y_0\not\in D$. Since $D=\mbox{Hlim} B(x_n,r_n)$, i.e.,
$\delta(\overline {B^*}\setminus B(x_n,r_n),\overline {B^*}\setminus
D)\rightarrow 0$. It follows from  Lemma 2.1 that there exists a
sequence $\{y_n\}$ satisfying
\begin{eqnarray}\label{2.11}
y_n\in \overline {B^*}\setminus B(x_n,r_n)\;\;\mbox{and}\;\;
y_n\rightarrow y_0.\
\end{eqnarray}
Hence
$$
d(y_n,x_n)\geq r_n\geq r_0.
$$
By (\ref{2.9}) and (\ref{2.11}) we may pass to the limit for
$n\rightarrow\infty$ to get
$$
d(y_0,x_0)\geq r_0,
$$
which leads to a contradiction and implies (\ref{2.10}) as desired.
This completes the proof. \hfill$\Box$

\vskip 0.1cm

{\bf Lemma 2.6}\quad Let $\{K_n\}_{n=1}^\infty$ be a sequence of
connected compact sets and $\mbox{H}\!\lim K_n=K$, where $K_n\subset
B^*$ for all $n\in \mathbb{Z}^+$, then $K$ is also connected compact
set.

{\bf Proof}\quad By Lemma 2.4, we obtain that $K$ is compact.

Now we show that $K$ is a connected set. Otherwise, there exist at
least two components in $K$, denote one of these by $K_1$, and
$K_2=K\setminus K_1$, then $K_1, K_2$ are compact, which shows that
$l_0\equiv\mbox{dist}(K_1, K_2)>0$. Since $\displaystyle K_n\subset
\bigcup_{z\in K}B(z, \frac{l_0}{4})={\big [}\bigcup_{z\in K_1}B(z,
\frac{l_0}{4}){\big ]}\cup {\big [}\bigcup_{z\in K_2}B(z,
\frac{l_0}{4}){\big ]}$ for $n$ large enough, following from
$\displaystyle {\big [}\bigcup_{z\in K_1}B(z, \frac{l_0}{4}){\big
]}\cap {\big [}\bigcup_{z\in K_2}B(z, \frac{l_0}{4}){\big
]}=\emptyset$ and connectedness of $K_n$ we obtain the
contradiction. \hfill$\Box$

\vskip 1mm

{\bf Lemma 2.7}\quad Let $\Omega$ be a bounded open set and
$d(K)\equiv\mbox{dist}(K, \partial\Omega)$, where $K\subset\Omega$
be compact sets. Then for all compact sets $K_1, K_2\subset\Omega$,
we have
$$|d(K_1)-d(K_2)|\leq\delta(K_1, K_2).$$

{\bf Proof}\quad We assume $d(K_1)\leq d(K_2)$, then we only need to
show $d(K_2)\leq d(K_1)+\delta(K_1, K_2)$. Since $K_1,
\partial\Omega$ are compact, there exist $x_1\in K_1, y_1\in
\partial\Omega$ such that $d(K_1)=|x_1-y_1|$. Note that
$\displaystyle\mbox{dist}(x_1, K_2)\leq \sup_{z\in
K_1}\mbox{dist}(z, K_2)\leq\delta(K_1, K_2)$, we obtain
$d(K_2)=\mbox{dist}(K_2, \partial\Omega)\leq \mbox{dist}(K_2,
y_1)\leq \mbox{dist}(K_2, x_1)+|x_1-y_1|\leq d(K_1)+\delta(K_1,
K_2)$. \hfill$\Box$

The main result in this paper is as follows.

\vskip 1mm

{\bf Theorem 2.1}\quad If $\{\Omega_m\}_{m=1}^\infty\subset
\mathscr{C}_{M,R}$, then there exists a subsequence
$\{\Omega_{m_k}\}_{k=1}^\infty$ of
$\displaystyle\{\Omega_m\}_{m=1}^\infty$ such that
$$\mbox{H}\!\!\!
\lim_{k\rightarrow \infty}\Omega_{m_k}=\Omega\ \ \mbox{and}\ \
\Omega\in \mathscr{C}_{M,R}.$$ i.e., $(\mathscr{C}_{M,R}, \rho)$ is
a compact metric space.

\vskip 1mm

{\bf Proof}\quad By Lemma 2.4 we know that there exists $\Omega\subset B^*$ and  a
subsequence
$\{\Omega_{m_k}\}_{k=1}^\infty$ of
$\displaystyle\{\Omega_m\}_{m=1}^\infty$such
that $\displaystyle \mbox{H}\!\!
\lim_{k\rightarrow \infty}\Omega_{m_k}=\Omega$, without loss of generality, we still denote
$\displaystyle \mbox{H}\!\!
\lim_{m\rightarrow \infty}\Omega_{m}=\Omega$. So we only need to prove that
$\Omega\in \mathscr{C}_{M,R}$.

In order to prove that $\Omega\in \mathscr{C}_{M,R}$, we separate it to the following steps.

$1^\circ.$\quad $\Omega$ is not empty.

By the definition of $\mathscr{C}_{M,R}$, there exists $B(x_m,
R)\subset \Omega_m$ for every $m\in \mathbb{Z}^+$. Now we consider
the sequence of open sets $\{B(x_m, R)\}_{n=1}^\infty$, by Lemma 2.4 and Lemma 2.5
we obtain that there exists $x_0\in \Omega$ such that $B(x_0,
R)\subset\mbox{Hlim}B(x_n,R)=D\subset \Omega$ and $x_{n}\rightarrow
x_0 (n\rightarrow \infty)$. Let $\tilde A_m=\overline{B^*}\setminus
B(x_m, R), \tilde A=\overline{B^*}\setminus D; A_m=\overline
{B^*}\setminus \Omega_m, A=\overline{B^*}\setminus \Omega$. Then
$\delta(A_m,A)\rightarrow 0$ and $\delta(\tilde A_m,\tilde
A)\rightarrow 0$, by Lemma 2.2 we get $B(x_0, R)\subset
D\subset\Omega$ and which show that $\Omega$ is not empty.

$2^\circ.$\quad $\Omega$ is connected.

Otherwise, there exists at least two components in $\Omega$, we
denote two of components of $\Omega$ by $\Omega_1, \Omega_2$.
Obviously, $\Omega_1\cap\Omega_2=\emptyset$. Since $\Omega_1,
\Omega_2$ are all open sets, there exist two points $x_1\in \Omega_1, x_2\in\Omega_2$ and
$d>0$ such that $B(x_1, d)\subset\Omega_1$ and $B(x_2,
d)\subset\Omega_2$, here $d^*=\min\{\mbox{dist}(x_1,
\partial\Omega_1), \mbox{dist}(x_2,
\partial\Omega_2)\}$. For any $0<r_0<d^*$, by Lemma 2.3 there exists $N_0>0$ such that for
all $m\geq N_0$ we have $\overline{B(x_1, r_0)\cup B(x_2,
r_0)}\subset \Omega_m$, furthermore, ${B(x_1, r_0)\cup B(x_2,
r_0)}\subset \Omega_m$. By the definition of $\mathscr{C}_{M,R}$,
there exists a connected compact set $K_m$ for every $m\geq N_0$
such that $x_1, x_2\in K_m$ and $\displaystyle\bigcup_{z\in K_m}B(z,
\frac{r_0}{M})\subset\Omega_m$. Since $K_m\not\subset
\Omega_1\cup\Omega_2$ by $\Omega_1\cap\Omega_2=\emptyset$, there
exists some $y_m\in K_m\setminus\Omega$ for every $m\geq N_0$, then
we obtain a sequence of points $\{y_m\}_{m=N_0}^\infty$ and
$\{y_m\}_{m=N_0}^\infty\cap\Omega=\emptyset$. Since
$\{y_m\}_{m=N_0}^\infty\subset B^*$ and $B^*$ is bounded, there
exists a subsequence of $\{y_m\}_{m=N_0}^\infty$, denote by itself,
and $y_0\not\in \Omega$ such that $y_m\rightarrow y_0$ as
$m\rightarrow \infty$.

Now we consider the ball $\displaystyle B(y_0, \frac{r_0}{2M})$.
Since $y_m\rightarrow y_0$ as $m\rightarrow \infty$, there exists
$N_1>0$ such that for all $m\geq N_1$ we have
$\displaystyle|y_m-y_0|<\frac{r_0}{2M}$. On the other hand, by
$\displaystyle\bigcup_{z\in K_m}B(z, \frac{r_0}{M})\subset\Omega_m$
we get $\displaystyle B(y_m, \frac{r_0}{M})\subset\Omega_m$,
hence, $\displaystyle B(y_0, \frac{r_0}{2M})\subset\Omega_m$
for all $m\geq N=\max\{N_0, N_1\}$. By the same argument with step
$1^\circ$ or by Lemma 2.2 we get $\displaystyle B(y_0,
\frac{r_0}{2M})\subset\Omega$, hence $\displaystyle B(y_0,
\frac{d^*}{2M})\subset\Omega$ by letting $r_0\rightarrow d^*$, which contradict to $y_0\not\in
\Omega$ and we prove that $\Omega$ is connected.

$3^\circ$.\quad $\Omega\in \mathscr{C}_{M,R}$.

We only need to prove that $\Omega$ satisfies the property ($C_M$).

(1). For any $x, y\in \Omega$, we set
\begin{eqnarray*}
\displaystyle d_0\equiv d_0(x,y)\equiv
\sup\{d(K); \  K\ \mbox{is connected compact set and}\ x, y\in K,\
\bigcup_{z\in K}B(z, d)\subset\Omega\}.
\end{eqnarray*}
We note that $\Omega$ is a local path connected set (since $\Omega$ is an open set)
and connected (by step $2^\circ$), then
$\Omega$ is a path connected set by the point topological theory. Hence
for any $x,y\in\Omega$ there exist at least one path $f: [0,1]\rightarrow \Omega$ such
that $f(0)=x, f(1)=y$. Which implies $d_0=d_0(x,y)$ exists.

Let
$\{K_m\}_{m=1}^\infty$ be a sequence such that $d(K_m)\rightarrow
d_0$ as $m\rightarrow\infty$ by Lemma 2.4. Since $K_m\subset B^*,\
m\in\mathbb{Z}^+$, there exist a set $K_0$ and a subsequence of
$\{K_m\}_{m=1}^\infty$, still denote by itself, such that $\delta(K_m,
K_0)\rightarrow 0$. By Lemma 2.7 we know that $d(K)$ is a continuous
function about $K$ under Hausdorff distance, hence
$d(K_0)=d_0$.

Following from Lemma 2.6 we get $K_0$ is a connected compact set,
and by Lemma 2.1 we obtain that $x, y\in K_0$.

(2). We claim that $\displaystyle \bigcup_{z\in K_0}B(z,
d_0)\subset\Omega$.

If the Claim is false, there exists a point $\displaystyle
x_0\in{\big [}\bigcup_{z\in K_0}B(z, d_0){\big ]}\setminus\Omega$.
We assume $\mbox{dist}(x_0, K_0)=|x_0-z_0|$, here $z_0\in K_0$, and
denote $\epsilon_0=d_0-|x_0-z_0|$, then $\displaystyle B(x_0,
\frac{\epsilon_0}{4})\subset\bigcup_{z\in K_0}B(z, d_0-\frac{3\epsilon_0}{4})$. But
$\displaystyle K_0\subset \bigcup_{z\in
K_m}B(z,\frac{\epsilon_0}{4})$ for $m$ large enough (see Remark 2.1), we get
$\displaystyle B(x_0, \frac{\epsilon_0}{4})\subset \bigcup_{z\in
K_m}B(z,d(K_m))\subset\Omega_m$ for $\displaystyle
d_0-d(K_m)<\frac{\epsilon_0}{2}$, then by the same argument with step
$1^\circ$ we obtain $\displaystyle B(x_0,
\frac{\epsilon_0}{4})\subset\Omega$. Which is a contradiction since
$x_0\not\in\Omega$.

(3). We show that $\displaystyle d_0\geq\frac{d^*}{M}$, where
$d^*=\min\{\mbox{dist}(x,\partial\Omega), \mbox{dist}(y,
\partial\Omega)\}$.

Since $\overline{B(x, r)}, \overline{B(y, r)}\subset\Omega$ for any
given $0<r<d^*$, there exists $N(r)>0$ such that for all $m\geq
N(r)$ we have $\overline{B(x, r)}, \overline{B(y,
r)}\subset\Omega_m$. Hence there exists $Q_m\subset\Omega_m$ such
that $\displaystyle \bigcup_{z\in
Q_m}B(z,\frac{r}{M})\subset\Omega_m$ for every $m\geq N(r)$ by
$\Omega_m\in\mathscr{C}_{M,R}$, here $Q_m$ is a connected compact
set. By Lemma 2.4 we obtain that there exist a set $Q_0$ and a
subsequence of $\{Q_m\}_{m=N(r)}^\infty$, still denote by itself, such
that $\delta(Q_m, Q_0)\rightarrow 0$ as $n\rightarrow\infty$. By
Lemma 2.6 we know that $Q_0$ is also a connected compact set.

Now claim: $\displaystyle \bigcup_{z\in
Q_0}B(z,\frac{r}{M})\subset\Omega$.

Otherwise, we take $\displaystyle x^*\in{\big [}\bigcup_{z\in
Q_0}B(z,\frac{r}{M}){\big]}\setminus\Omega$, then there exists
$z^*\in Q_0$ such that $|x^*-z^*|=\mbox{dist}(x, Q_0)$. Denote
$\displaystyle\epsilon^*=\frac{r}{M}-|x^*-z^*|$, then $\displaystyle B(x^*,
\frac{\epsilon^*}{2})\subset \bigcup_{z\in
Q_0}B(z,\frac{r}{M}-\frac{\epsilon^*}{2})$. By Remark 2.1 we know
that $\displaystyle Q_0\subset\bigcup_{z\in
Q_m}B(z,\frac{\epsilon^*}{2})$ for all $m\geq N^*=\max\{N(r),
N(\epsilon^*)\}$, and we get $\displaystyle B(x^*,
\frac{\epsilon^*}{2})\subset \bigcup_{z\in
Q_m}B(z,\frac{r}{M})\subset\Omega_m$ for all $m\geq N^*$, by the
same argument with step $1^\circ$ we obtain $\displaystyle B(x^*,
\frac{\epsilon^*}{2})\subset\Omega$. Which is a contradiction since
$x^*\not\in\Omega$.

Since $0<r<d^*$ is arbitrary, letting $r\rightarrow d^*$ we get
$\displaystyle \bigcup_{z\in Q_0}B(z,\frac{d^*}{M})\subset\Omega$.
By the definition of $d_0$ we obtain $\displaystyle
d_0\geq\frac{d^*}{M}$. Which implies $\Omega\in \mathscr{C}_{M,R}$.

$4^\circ$. Following from $1^\circ, 2^\circ, 3^\circ$, we have
proved the Theorem 2.1.\hfill$\Box$

\vskip 1mm

\section {Existence of shape optimization}

\quad\quad In this section, we shall prove the existence of
problems $(P)$.

\vskip 0.1cm

{\bf Theorem 3.1.}\quad The shape optimization problem $(P)$ has
at least one solution.

\vskip 0.1cm

{\bf Proof.}$\;\;\;$Throughout the proof of Theorem 3.1, we shall
use $spt(u)$ to denote the support of $u$.

Let $\displaystyle d=\inf_{\Omega\in \mathscr{C}_{M,R}} J(\Omega)=\inf_{\Omega\in\mathscr{C}_{M,R}}{\Big\{} \frac{1}{2}\int_{B^*}|u_\Omega-g|^2dx{\Big \}}$. It is obvious that $d>-\infty$. Then
there exists a sequence $\{\Omega_m\}_{m=1}^\infty\subset \mathscr{C}_{M,R}$
such that
\begin{eqnarray}\label{3.1}
d=\lim_{m \rightarrow \infty} \frac{1}{2}\int_{B^*}|u_m-g|^2dx ,
\end{eqnarray}
where $u_m\equiv u_{\Omega_m}$ is the  weak solution of (\ref{2}). By
Theorem 2.1, there exist a subsequence of
$\{\Omega_n\}_{n=1}^\infty$, still denoted by itself,
 and $\Omega^*\in \mathscr{C}_{M,R}$  such that
$\Omega^*=\mbox{Hlim}\Omega_n$.

By taking $u=u_m$, $\Omega=\Omega_m$ in (\ref{3}), we get
\begin{eqnarray*}
\int_{\Omega_m} \langle A.\nabla u_m,\nabla\phi\rangle dx
=\langle f|_{\Omega_m},\phi\rangle_{H^{-1}({\Omega_m})
\times H_0^1({\Omega_m})}\quad \forall \phi\in C_0^\infty({\Omega_m})
\end{eqnarray*}
Since $u_m\in H_0^1(\Omega_m)$, we have
\begin{eqnarray}\label{4.0}
\alpha\int_{\Omega_m}|\nabla u_m|^2dx&\leq&\int_{\Omega_m} \langle A.\nabla u_m,\nabla u_m\rangle\,dx \cr
&=&\langle f|_{\Omega_m},u_m\rangle_{H^{-1}({\Omega_m})
\times H_0^1({\Omega_m})}\cr
&\leq& \|f\|_{H^{-1}(B^*)}\cdot \|u_m\|_{H_0^1(\Omega_m)}
\end{eqnarray}
by (\ref{0}), which implies that
$$
\int_{\Omega_m} |\nabla u_m|^2\,dx\leq C,
$$
here and throughout the proof of Theorem 3.1, $C$ denotes several
positive constants independent of $m$. Let
\begin{eqnarray}\label{3.3}
{\hat{u}}_m(x)=\left\{
\begin{array}{ll}
u_m(x)&\mbox{in}\;\;\Omega_m,\\[3mm]
0&\mbox{in}\;\; B^*-\Omega_m,\\
\end{array}
\right.
\end{eqnarray}
then $\{\hat{u}_m\}_{m=1}^\infty$ is bounded in
$H_0^1(B^*)$. Hence there exists subsequence of
$\{\hat{u}_m\}$, still denoted by itself, such that
\begin{eqnarray}\label{3.4}
\hat{u}_m\rightarrow \hat{u}\;\mbox{weakly
in}\;H_0^1(B^*)\; \mbox{and strongly in}\; L^2(B^*)
\end{eqnarray}
for some $\hat{u}\in H_0^1(B^*)$.

Now we claim that
\begin{eqnarray}\label{3.5}
\hat{u}(x)\in H_0^1(\Omega^*).
\end{eqnarray}

In fact, we only need to show  that
\begin{eqnarray}\label{3.6}
\hat{u}(x)=0\;\mbox{a.e. in}\; B^*\setminus\overline{\Omega^*}.
\end{eqnarray}
Indeed, for any open subset $K$ satisfying $\overline K\subset B^*-\overline{\Omega^*}$,
denote $l_d=\mbox{dist}(\overline K,\overline{\Omega^*})$, since $\Omega^*=\mbox{Hlim}\Omega_n$, i.e.,
$\delta(\overline{B^*}\setminus \Omega_m, \overline{B^*}\setminus \Omega^*)\rightarrow 0$,
which implies that for every $\d 0<\frac{l_d}{4}$ there exists an integer $m_0>0$ such that
for any $n\geq n_0$ we have
$\displaystyle\delta(\overline{B^*}\setminus \Omega_m, \overline{B^*}\setminus \Omega^*)<\frac{l_d}{4}$.
Hence by the definition of $\delta(\cdot,\cdot)$ we have
$$\sup_{x\in \overline{B^*}\setminus \Omega}d(x, \overline{B^*}\setminus \Omega_m)
<\frac{l_d}{4},\ \forall \ m\geq m_0.$$
i.e., $\overline{B^*}\setminus \Omega^*\subset
\{y\in\mathbb{R}^k; \ d(y, \overline{B^*}\setminus \Omega_m< {\frac{l_d}{4}})\}$.
So we obtain
$$\{y\in \Omega_m;\ d(y,\partial\Omega_m)>
{\frac{l_d}{4}}\}\subset\Omega^*,
\ \forall \ m\geq m_0.$$
Furthermore, we have
$$\{y\in \mathbb{R}^k;\ d(y,\partial\Omega^*)<
{\frac{l_d}{2}}\},\ \forall \ m\geq m_0.$$
Hence $\overline K\subset B^*\setminus \overline{\Omega_m}, \ \forall \ m\geq m_0$.
Thus
$$
\int_K |\hat{u}(x)|^2\,dx=\lim_{m\rightarrow \infty} \int_K
|\hat{u}_m(x)|^2\,dx\leq \overline {\lim_{m\rightarrow \infty}}
\int_{B^*\setminus\overline{\Omega_m}} |\hat{u}_m(x)|^2\,dx=0,
$$
which implies that $\hat{u}(x)=0$ a.e. in $K$. Since
$K\subset \overline K\subset B^*\setminus\overline{\Omega^*}$ is arbitrary, (\ref{3.6})
and consequently
(\ref{3.5}) follow.

Then we claim that
\begin{eqnarray}\label{4.1}
\int_{\Omega^*} \langle A.\nabla \hat{u},\nabla\phi\rangle dx
=\langle f|_{\Omega^*},\phi\rangle_{H^{-1}({\Omega^*})\times H_0^1({\Omega^*})}\quad
\forall \phi\in C_0^\infty({\Omega^*})
\end{eqnarray}
i.e.,
\begin{eqnarray*}
\int_{spt(\phi)} \langle A.\nabla \hat{u},\nabla\phi\rangle dx
=\langle f|_{spt(\phi)},\phi\rangle_{H^{-1}({\Omega^*})\times H_0^1({\Omega^*})}\quad
\end{eqnarray*}
for each $\phi\in C_0^\infty({\Omega^*})$.

 Let
\begin{eqnarray}\label{3.8}
\hat{\phi}=\left\{
\begin{array}{ll}
\phi&\mbox{in}\;\;\Omega^*,\\[3mm]
0&\mbox{in}\;\;B^*-\overline {\Omega^*}, \\
\end{array}
\right.
\end{eqnarray}
by Lemma 2.3, there exists a positive integer $m_1(\phi)$,
such that
$$
spt(\hat{\phi})\;=spt(\phi)\subset
\Omega_m,\;\;\;\mbox{for all}\;\;m\geq m_1(\phi).
$$
Then for each $m\geq m_1(\phi),\  \hat{\phi}\in
C_{0}^\infty (\Omega_m)$. So by (\ref{3}), we have
\begin{eqnarray*}
\int_{\Omega_m} \langle A.\nabla \hat{u}_m,\nabla {\hat\phi}\rangle\,dx
=\langle f|_{\Omega_m},\hat{\phi}\rangle_{H^{-1}({\Omega_m})
\times H_0^1({\Omega_m})}
\end{eqnarray*}
which, together with (\ref{3.8}), implies
\begin{eqnarray*}
\int_{\Omega_m} \langle A.\nabla \hat{u}_m,\nabla {\phi}\rangle\,dx
=\langle f|_{\Omega_m},{\phi}\rangle_{H^{-1}({\Omega_m})
\times H_0^1({\Omega_m})}
\end{eqnarray*}
passing to the limit for $m\rightarrow \infty$ and using (\ref{3.4}) we
get (\ref{4.1}).

Finally, we claim that
\begin{eqnarray}\label{3.9}
d\geq \frac{1}{2}\int_{B^*}|\hat{u}-g|^2dx.
\end{eqnarray}
We notice that
\begin{eqnarray*}
\hat{u}_m\rightarrow \hat{u} \quad \mbox{strongly in}\quad L^2(B^*).
\end{eqnarray*}
Hence we have
\begin{eqnarray*}
\frac{1}{2}\int_{B^*}|\hat{u}_m-g|^2\rightarrow \frac{1}{2}\int_{B^*}|\hat{u}-g|^2.
\end{eqnarray*}
and by (\ref{3.1}) we get
\begin{eqnarray*}
d=\frac{1}{2}\int_{B^*}|\hat{u}-g|^2.
\end{eqnarray*}
i.e. (\ref{3.9}) holds.

Following from (\ref{3.5}), (\ref{4.1}) and (\ref{3.9}) we obtain that $\Omega^*$ is a solution of
problem $(P)$. This completes the proof. \hfill$\Box$

\end{document}